\documentclass[preprint,3p,12pt]{elsarticle}

\usepackage{lineno,hyperref}
\modulolinenumbers[5]

\usepackage{hyperref}

\usepackage{epstopdf}

\usepackage{mathrsfs}
\usepackage{amssymb}
\usepackage{amsfonts}
\usepackage{latexsym}
\usepackage{amsmath}
\usepackage{xcolor}
\usepackage{wrapfig}
\usepackage{floatflt}

\usepackage{graphicx}
\def\lb{\label}

\newcommand{\er}[1]{\textrm{(\ref{#1})}}

\newtheorem{theorem}{\bf Theorem}[section]

\newtheorem{definition}[theorem]{\bf Definition}

  \def\cA{{\mathcal A}}       \def\mA{{\mathscr A}}
          \def\mB{{\mathscr B}}

\def\d{\delta}         
\def\D{\Delta}         
           
   \def\cH{{\mathcal H}}       \def\mH{{\mathscr H}}
    \def\cI{{\mathcal I}}       
           
 \def\cK{{\mathcal K}}       \def\mK{{\mathscr K}}
         
\def\l{\lambda}

\def\s{\sigma}  \def\cR{{\mathcal R}}

\def\o{\omega}

\def\O{\Omega}

\def\Z{{\mathbb Z}}       \def\C{{\mathbb C}}    
    \def\N{{\mathbb N}}   

               \def\wt{\widetilde}
\def\no{\noindent}

\let\ge\geqslant                 \let\le\leqslant

\def\ss{\subset}                 \def\ts{\times}

\def\el2{\ell^{\,2}}             \def\1{1\!\!1}

\def\det{\mathop{\mathrm{det}}\nolimits}

\def\dim{\mathop{\mathrm{dim}}\nolimits}

\def\BBox{\hspace{1mm}\vrule height6pt width5.5pt depth0pt \hspace{6pt}}

\let\ge\geqslant
\let\le\leqslant

\newcommand{\ca}{\begin{cases}}
\newcommand{\ac}{\end{cases}}
\newcommand{\ma}{\begin{pmatrix}}
\newcommand{\am}{\end{pmatrix}}
\def\eq{\begin{equation}}
\def\qe{\end{equation}}
\def\[{\begin{equation}}
\def\]{\end{equation}}

\def\BBox{\hspace{1mm}\vrule height6pt width5.5pt depth0pt \hspace{6pt}}

\begin{document}

\begin{frontmatter}

\title{Algebra of 2D periodic operators with local and perpendicular defects}
\date{\today}

\author
{Anton A. Kutsenko}

\address{Department of Mathematics, Aarhus University, Aarhus, DK-8000,
Denmark; email: akucenko@gmail.com}

\begin{abstract}
We show that 2D periodic operators with local and perpendicular
defects form an algebra. We provide an algorithm of finding spectrum
for such operators. While the continuous spectral components can be
computed by simple algebraic operations on some matrix-valued
functions and few number of integrations, the discrete part is much
more complicated.
\end{abstract}

\begin{keyword}
periodic lattice with defects, Floquet - Bloch spectrum, guided
waves, localised waves (states), operator algebras
\end{keyword}

\end{frontmatter}


\section{Introduction}

Defects in periodic structures play a major role in various fields
of science, see, e.g., discussions in \cite{Kjmaa}. The algebras of
multidimensional discrete periodic operators with parallel defects
(APD) are considered in the mentioned paper, where the algorithm of
finding spectrum based on simple algebraic operations and few number
of integration is provided (we call such algorithms as explicit). In
the current article we extend two-dimensional APD by adding
perpendicular defects. Generally speaking, in comparison with the
parallel defects the perpendicular defect makes non-explicit the
algorithm of finding point spectrum. At the same time, the
continuous components can be computed explicitly in the same way as
in \cite{Kjmaa}, \cite{K1a}. As well as parallel defects, the
perpendicular defects have a lot of applications, see, e.g.,
\cite{WHMH}, \cite{CWC}, \cite{JMFVJH}, \cite{JJ} about electro and
optical crossing wave-guides. Some comparison of parallel and
perpendicular waveguides is treated in \cite{HBCOASLXWMH}. The
methods of finding guided and local waves, and the corresponding
spectrum are usually approximative and are based on supercell
approaches, where the infinite structure is replaced with a large
finite structure which has a discrete spectrum only. In the current
paper we propose a non-approximative algorithm of finding spectrum
based on an expansion of the periodic operator with defects into the
product of the operators with "simple" spectral components.

\begin{figure}[h]
\begin{minipage}[h]{0.49\linewidth}
\center{\includegraphics[width=0.99\linewidth]{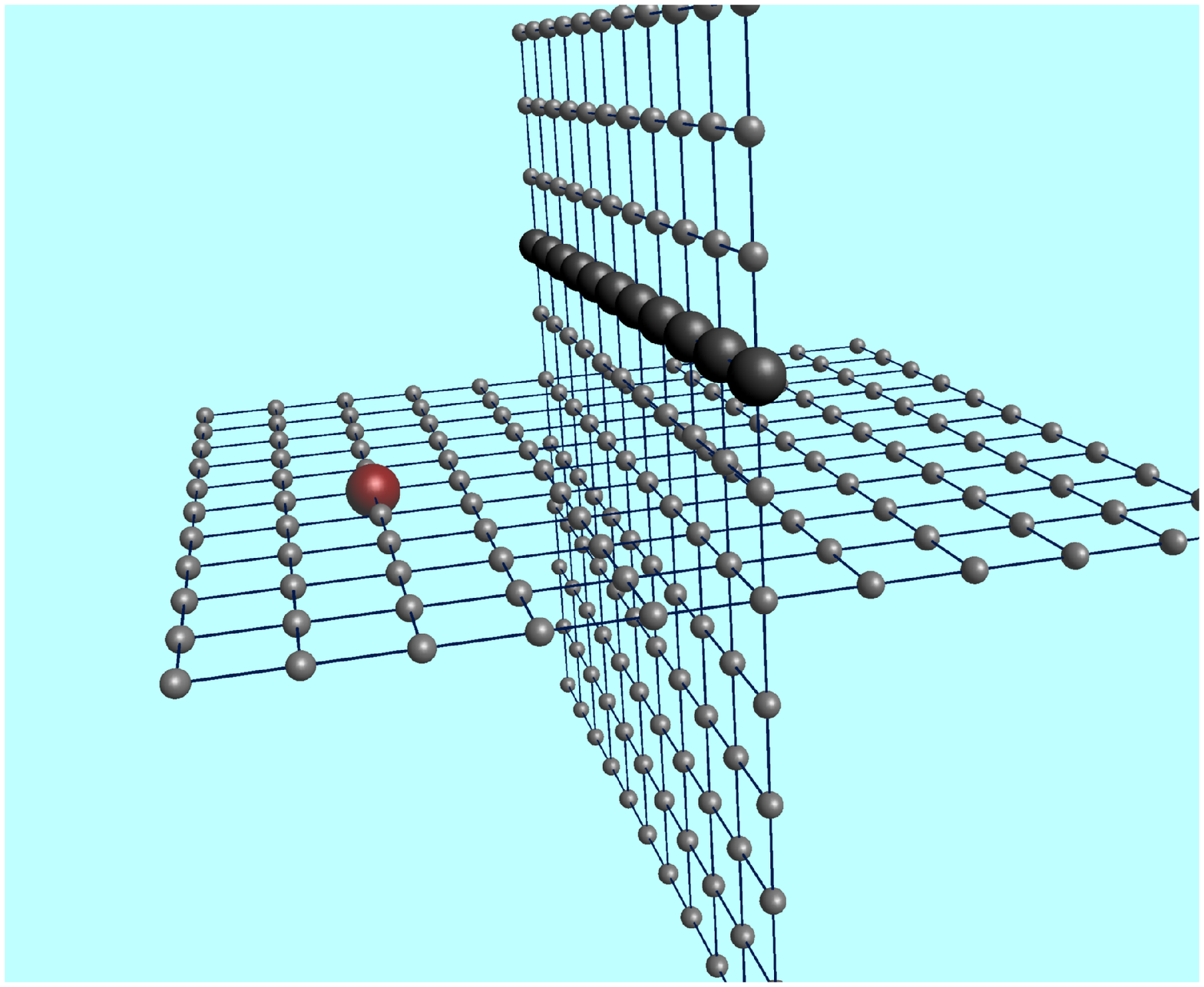}} \\ (a)
\end{minipage}
\hfill
\begin{minipage}[h]{0.49\linewidth}
\center{\includegraphics[width=0.99\linewidth]{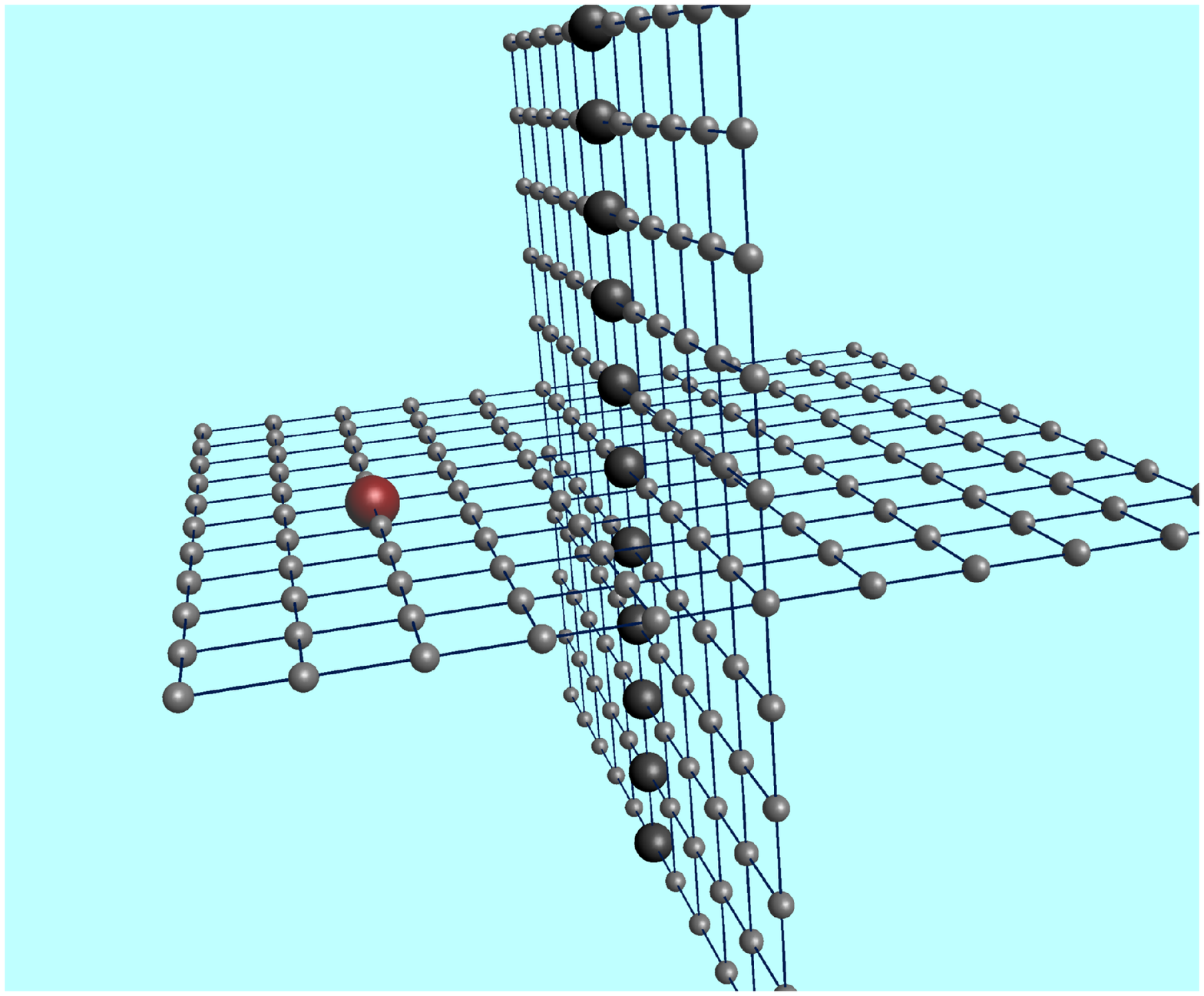}} \\ (b)
\end{minipage}
\caption{ Two coupled 2D lattices with parallel (a) and
perpendicular (b) line defects, and one point defect. Periodic
operators on the first structure belong to $\mH_{2,M}$, on the
second belong to $\mA_{2,M}$.} \label{fig1}
\end{figure}

Let $M$ be some positive integer. Introduce the following Hilbert
space and integral operators
\[\lb{001}
 L^2_{2,M}:=L^2([0,1]^2,\C^M),\ \ \langle\cdot\rangle_i:=\int_0^1\cdot
 dk_i,\ i=1,2,\ \
 \langle\cdot\rangle_{12}:=\langle\langle\cdot\rangle_1\rangle_2,
\]
where $\cdot$ means some matrix- or vector-valued function depending
on ${\bf k}=(k_1,k_2)\in[0,1]^2$. The bounded and compact operators
acting on $L^2_{2,M}$ will be denoted as $\mB_{2,M}$ and $\mK_{2,M}$
respectively. By analogy with APD (see \cite{Kjmaa}) define

\begin{definition}\lb{D1} The algebra of 2D periodic operators with local and perpendicular defects
\[\lb{002}
 \mA_{2,M}={\rm Alg}\{{\bf A}\cdot,\langle\cdot\rangle_1,\langle\cdot\rangle_2,\cK\}
\]
is a minimal non-closed subalgebra of the algebra $\mB_{2,M}$ which
contains all operators of multiplication by $M\ts M$ continuous
matrix-valued functions ${\bf A}\cdot$, the integral projectors
$\langle\cdot\rangle_i$, $i=1,2$, and all compact operators
$\cK\in\mK_{2,M}$. Here $\cdot$ denotes the operator argument ${\bf
u}\in L^2_{2,M}$.
\end{definition}

Note that APD in 2D case have the form
\[\lb{003}
 \mH_{2,M}={\rm Alg}\{{\bf
 A}\cdot,\langle\cdot\rangle_1,\langle\cdot\rangle_{12}\},\ \ \wt\mH_{2,M}={\rm Alg}\{{\bf
 A}\cdot,\langle\cdot\rangle_2,\langle\cdot\rangle_{12}\}.
\]
The general class of multidimensional APD is studied in \cite{Kjmaa}
and \cite{K1a}. The algebra \er{002} is an extension of APDs
\er{003}, i.e.
\[\lb{004}
 \mH_{2,M}\ss\mA_{2,M},\ \ \wt\mH_{2,M}\ss\mA_{2,M}.
\]
As mentioned above, this extension makes non-explicit the algorithm
of finding eigenvalues. The schematic difference between some
operators from $\mH_{2,M}$ and $\mA_{2,M}$ is illustrated in Fig.
\ref{fig1}. The next theorem is some analogue of the theorem from
\cite{Kjmaa}.

\begin{theorem}\lb{P1}
Each operator $\cA\in\mA_{2,M}$ has a following representation
\[\lb{005}
 \cA {\bf u}={\bf A}_0{\bf u}+{\bf A}_1\langle{\bf B}_1{\bf
 u}\rangle_{1}+{\bf A}_2\langle{\bf B}_2{\bf
 u}\rangle_{2}+\cK{\bf u},\ \ {\bf
 u}\in L^2_{2,M},
\]
where $\cK\in\mK_{2,M}$ and ${\bf A}$, ${\bf B}$ are continuous
matrix-valued functions on $[0,1]^2$ of sizes
\[\lb{006}
 \dim({\bf A}_0)=M\ts M,\ \ \dim({\bf B}_j)=M_j\ts M,\ \ \dim({\bf
 A}_j)=M\ts M_j,\ \ j=1,2
\]
with some positive integers $M_j$. The set of all operators of the
form \er{002} coincides with $\mA_{2,M}$.
\end{theorem}

To describe the spectrum of the operators \er{005} we need the
following definition.

\begin{definition}\lb{D2}
Let $\cA\in\mA_{2,M}$ be some operator of the form \er{002}. Define
the following matrix-valued functions (if the inverse matrices
exist)
\[\lb{007}
 {\bf C}_0={\bf A}_0^{-1},\ \ {\bf E}_0={\bf A}_0,\ \ {\bf C}_1={\bf C}_0{\bf A}_1,\ \
 {\bf C}_2={\bf C}_0{\bf A}_2,
\]
\[\lb{008}
 {\bf E}_1={\bf I}+\langle{\bf B}_1{\bf
 C}_1\rangle_{1},\ \ {\bf E}_2={\bf I}+\langle{\bf B}_2{\bf
 C}_2\rangle_{2},
\]
\[\lb{009}
 {\bf D}_1({\bf k},{\bf k}')={\bf C}_1({\bf k}){\bf E}_1^{-1}(k_2){\bf B}_1(k_1',k_2){\bf C}_2(k_1',k_2){\bf
 B}_2({\bf k}'),
\]
\[\lb{010}
 {\bf D}_2({\bf k},{\bf k}')={\bf C}_2({\bf
 k}){\bf E}_2^{-1}(k_1)\int_0^1{\bf B}_2(k_1,k_2''){\bf D}_1(k_1,k_2'',{\bf
 k}')dk_2'',
\]
where ${\bf k}=(k_1,k_2)$ and ${\bf k}'=(k_1',k_2')$.
\end{definition}

The next Theorem is our main result. It not only provides the
explicit procedure of finding inverse operators, but along with the
Theorem \ref{P1} shows that the subset of $\mA_{2,M}$ consisting of
all invertible operators $({\rm Inv}\mA_{2,M},\circ)$ is an
algebraic group with the multiplication given by the composition of
mappings (the multiplication is the same as in $\mA_{2,M}$).
Theoretically, it could be that $\cA^{-1}\not\in\mA_{2,M}$ for
$\cA\in\mA_{2,M}$, but it did not happen. Everywhere "inverse" means
inverse in the large algebra $\mB_{2,M}$ of all bounded operators.

\begin{theorem}\lb{T1}
An operator $\cA\in\mA_{2,M}$ of the form \er{005} is invertible if
and only if $\det{\bf E}_j\ne0$ \er{007}-\er{008} everywhere for
$j=0,1,2$ and the operator $\cI+\cK_1$ is invertible. The operator
$\cI$ is the identity operator and the compact operator $\cK_1$ is
defined by
\[\lb{011}
 \cK_1{\bf u}=\int_{[0,1]^2}({\bf D}_2-{\bf D}_1)({\bf k},{\bf k}'){\bf
 u}({\bf k}')d{\bf k}'+\cR\circ\cK{\bf u},\ \ {\bf u}\in L^2_{2,M},
\]
\[\lb{012}
 \cR=(\cI-{\bf C}_2{\bf E}_2^{-1}\langle{\bf
 B}_2\cdot\rangle_2)\circ(\cI-{\bf C}_1{\bf E}_1^{-1}\langle{\bf
 B}_1\cdot\rangle_1)\circ({\bf C}_0\cdot),
\]
where $\cdot$ means operator argument. Moreover, the inverse
operator $\cA^{-1}\in\mA_{2,M}$ has the form
\[\lb{inva}
 \cA^{-1}=(\cI-\cK_1\circ(\cI+\cK_1)^{-1})\circ\cR.
\]
\end{theorem}

The Theorem \ref{T1} immediately yields to the next Corollary
describing the spectrum.

\no{\bf Corollary.} {\it Let $\cA\in\mA_{2,M}$ be an operator of the
form \er{005}. Taking ${\bf A}_0:={\bf A}_0-\l{\bf I}$ (${\bf I}$ is
the identity matrix) in the procedure \er{007}-\er{010} we obtain
that the spectrum of $\cA$ is
\[\lb{013}
 \s(\cA)=\bigcup_{j=0}^3\s_j\ \ with\ \ \s_0=\{\l:\ \det{\bf E}_0=0\ for\ some\ {\bf
 k}\in[0,1]^2\},
\]
\[\lb{014}
 \ca \s_1=\{\l:\ \det{\bf E}_1(k_2)=0\ fro\ some\ k_2\in[0,1]\},\\ \s_2=\{\l:\ \det{\bf E}_2(k_1)=0\ fro\ some\
 k_1\in[0,1]\},\ac
\]
\[\lb{015}
 \s_3=\{\l:\ \cI+\cK_1\ is\ non-invertible\}.
\]
}

\no {\bf Remark.} 1) The matrix-valued function ${\bf E}_0$ is
defined for any $\l\in\C$. The matrix-valued functions ${\bf E}_1$,
${\bf E}_2$ are well-defined for $\l\not\in\s_0$. More precisely, we
can define ${\bf E}_{1,2}$ for some $\l\in\s_0$ but it does not
affect the spectrum as a set. The analytic (compact-)operator-valued
function $\cK_1(\l)$ is well-defined for
$\l\not\in\s_0\cup\s_1\cup\s_2$. So, the procedure of finding
spectrum consists of determining $\s_0$, then $\s_1$, $\s_2$, and
then $\s_3$.

2) As for the parallel defects (see \cite{Kjmaa}), $\s_0$
corresponds to non-attenuated eigensolutions, $\s_1,\s_2$ correspond
to guided eigensolutions, and $\s_3$ are eigenvalues. Note that
$\s_0$ does not depend on any perturbation of lower dimension,
$\s_1$ and $\s_2$ do not dependent on each other and on the compact
perturbation $\cK$.

3) In Corollary, instead of ${\bf A}_0:={\bf A}_0-\l{\bf I}$ we may
assume a general situation of extended spectral problems where all
${\bf A}_j$, ${\bf B}_j$  somehow depend on the spectral parameter
$\l$.

4) Along with \er{inva} we have the decomposition of the direct
operator
\[\lb{016}
 \cA=({\bf A}_0\cdot)\circ(\cI+{\bf C}_2\langle{\bf B}_2\cdot\rangle_{2})\circ
 (\cI+{\bf C}_1\langle{\bf B}_1\cdot\rangle_{1})\circ(\cI+\cK_1).
\]
Each term belongs to the corresponding subgroup of invertible
operators from $\mA_{2,M}$. This decomposition is unique. The
situation is similar to that of \cite{K1a} except that we probably
can not correctly define the vector-valued traces and determinants
of $\cA$.


The work is organized as follows: Section \ref{S1} contains proofs
of our results. Section \ref{S2} provides an application of our
results to the problem of wave propagation through 2D spring-mass
model with two perpendicular wave-guides. The conclusion is given in
Section \ref{S3}.

\section{\lb{S1} Proof of Theorems \ref{P1}, \ref{T1}}

{\bf Proof of Theorem \ref{P1}.} It is obvious that any $\cA$
\er{005} belongs to $\mA_{2,M}$. To complete the proof we need to
show that the sum and products of the operators of the form \er{005}
have the same form. It is sufficient to show this fact for summands
only. For the components from $\cH_{2,M}$ and $\wt\cH_{2,M}$ the
corresponding identities  are already shown in \cite{Kjmaa}. It is
also obvious that $\mK_{2,M}$ is a two-sided ideal in $\mA_{2,M}$.
It remains to show that:
\[\lb{100}
 ({\bf A}_1\langle{\bf B}_1\cdot\rangle_{1})\circ({\bf A}_2\langle{\bf
 B}_2{\bf u}\rangle_{2})=\int_{[0,1]^2}{\bf D}({\bf k},{\bf k}'){\bf
 u}({\bf k}')d{\bf k}'
\]
is obviously a compact operator with the continuous kernel
\[\lb{101}
 {\bf D}({\bf k},{\bf k}')={\bf A}_1({\bf k}){\bf B}_1(k_1',k_2){\bf A}_2(k_1',k_2){\bf
 B}_2({\bf k}').
\]
If we change the multipliers in \er{100} then we also obtain a
compact operator. \BBox

{\bf Proof of Theorem \ref{T1}.} Suppose that ${\bf A}_0$ is
non-invertible for some ${\bf k}^0\in[0,1]^2$ with a corresponding
null-vector ${\bf f}\in\C^M$ having the unit Euclidean norm.
Consider some sequence of characteristic functions $\chi_n({\bf k})$
of the sets $\O_n\ss[0,1]^2$, $n\in\N$, where domains $\O_n$ tends
to the point ${\bf k}^0$. Taking constants $c_n=({\rm
mes}\O_n)^{-\frac12}$ we obtain that ${\bf u}_n=c_n\chi_n{\bf f}$
have unit $L^2_{2,M}$-norm and
\[\lb{102}
 {\bf A}_0{\bf u}_n\to0,
\]
since ${\bf f}$ is a null-vector of the continuous matrix-valued
function ${\bf A}_0$ at the point ${\bf k}^0$. At the same time, in
\cite{Kjmaa} it is proved that $\O_n$ can be chosen such that
\[\lb{103}
 {\bf A}_i\langle{\bf B}_i{\bf u}_n\rangle_i\to0,\ \ i=1,2.
\]
Let ${\bf u}\in L^2_{2,M}$ be some continuous function. Then the
$L^2_{2,M}$-inner product ($^*$ means Hermitian conjugation)
\[\lb{104}
 \langle{\bf u}^*{\bf u}_n\rangle_{12}=({\rm mes}\O_n)^{\frac12}{\bf u}({\bf k}^0)^*{\bf
 f}+o(1)\to0
\]
since $\O_n$ tends to the point ${\bf k}^0$ and hence their Lebesgue
measures tends to 0. Then $\langle{\bf u}^*{\bf
u}_n\rangle_{12}\to0$ for any ${\bf u}\in L^2_{2,M}$ since
$L^2_{2,M}$-norm  of ${\bf u}_n$ is bounded (equal to $1$). Since
$\cK$ is compact we may assume that $\cK{\bf u}_n\to{\bf v}$, ${\bf
v}\in L^2_{2,M}$. Then
\[\lb{105}
 0=\lim\langle(\cK^*{\bf v})^*{\bf u}_n\rangle_{12}=\lim\langle{\bf v}^*\cK{\bf u}_n\rangle_{12}=\langle{\bf v}^*{\bf v}\rangle_{12}
\]
or in other words
\[\lb{106}
 \cK{\bf u}_n\to0.
\]
Identities \er{102}, \er{013}, and \er{106} leads to
$
 \cA{\bf u}_n\to0
$ which with $\|{\bf u}_n\|=1$ ($\|\cdot\|$ is $L^2_{2,M}$ norm)
means that $\cA$ is non-invertible (by Banach theorem about bounded
inverse linear mappings).

Suppose that ${\bf E}_0$ is invertible everywhere. Then $\cA$ and
\[\lb{107}
 \cA_1=({\bf E}_0\cdot)^{-1}\circ\cA=({\bf C}_0\cdot)\circ\cA=\cI+{\bf C}_1\langle{\bf B}_1\cdot\rangle_{1}+
 {\bf C}_2\langle{\bf B}_2\cdot\rangle_{2}+({\bf C}_0\cdot)\circ\cK
\]
are invertible or non-invertible simultaneously. If ${\bf E}_1$ is
non-invertible at some $k_2^0\in[0,1]$ then as it is shown in
\cite{Kjmaa} the operator $\cI+{\bf C}_1\langle{\bf
B}_1\cdot\rangle_{1}$ is non-invertible and there exist domains
$\O_n\ss[0,1]$ tending to $k_2^0$  such that
\[\lb{108}
 (\cI+{\bf C}_1\langle{\bf B}_1\cdot\rangle_{1}){\bf u}_n\to0,
\]
where ${\bf u}_n=c_n\chi_n{\bf C}_1{\bf f}$, $\chi_n({\bf k})$ is
the characteristic function of the set $[0,1]\ts\O_n$, ${\bf f}$ is
a null-vector of ${\bf E}_1(k_2^0)$ with the unite Euclidean norm,
and $c_n$ are taken such that $\|{\bf u}_n\|=1$. It is true that for
some $\wt k_1$ the Euclidean norm of ${\bf C}_1(\wt k_1,k_2^0){\bf
f}$ is non-zero since otherwise
\[\lb{109}
 {\bf 0}={\bf E}_1(k_2^0){\bf f}={\bf f}.
\]
Then we have
\[\lb{110}
 1=\|{\bf u}_n\|\ge\d({\rm mes}\O_n)^{\frac12}c_n,
\]
where the absolute constant $\d$ depends on the the Euclidean norm
of ${\bf C}_1(\wt k_1,k_2^0){\bf f}$ only (recall that ${\bf C}_1$
is a continuous matrix-valued function). The following estimates are
fulfilled
\[\lb{111}
 \|{\bf C}_2\langle{\bf B}_2{\bf u}_n\rangle_{2}\|\le C({\rm
 mes}\O_n)c_n\le (C/\d)({\rm mes}\O_n)^{\frac12},
\]
where $C$ is an absolute constant depending on ${\bf C}_1$, ${\bf
C}_2$, and ${\bf B}_2$. The right-hand side of \er{111} tend to $0$
since $\O_n$ tend to the point $k_2^0$. Using this fact, \er{108},
and $({\bf C}_0\cdot)\circ\cK{\bf u}_n\to0$ (which can be proved in
the same manner as \er{106}) we deduce that the operator $\cA_1$ is
non-invertible by Banach theorem about bounded inverse linear
mappings.

Suppose that ${\bf E}_1$ is invertible everywhere. Then it is not
difficult to show that (see \cite{Kjmaa})
\[\lb{112}
 (\cI+{\bf C}_1\langle{\bf B}_1\cdot\rangle_{1})^{-1}=\cI-{\bf C}_1{\bf E}_1^{-1}\langle{\bf
 B}_1\cdot\rangle_{1}.
\]
Multiplying \er{107} by \er{112} we deduce that the operators
$\cA_1$ and
\[\lb{113}
 \cA_2=\cI+{\bf C}_2\langle{\bf B}_2\cdot\rangle_{2}-\int_{[0,1]^2}{\bf D}_1({\bf k},{\bf k}')\cdot({\bf k}')d{\bf
 k}'+(\cI-{\bf C}_1{\bf E}_1^{-1}\langle{\bf
 B}_1\cdot\rangle_{1})\circ({\bf C}_0\cdot)\circ\cK
\]
are invertible or non-invertible simultaneously. Here $\cdot$ means
an operator argument (as usual). Now we can apply again the above
arguments. If ${\bf E}_2$ is non-invertible at some $k_1^0\in[0,1]$
then the operator $\cI+{\bf C}_2\langle{\bf B}_2\cdot\rangle_{2}$ is
non-invertible and hence $\cA_2$ is non-invertible. Suppose that
${\bf E}_2$ is invertible everywhere. Then $\cI+{\bf C}_2\langle{\bf
B}_2\cdot\rangle_{2}$ is invertible with
\[\lb{114}
 (\cI+{\bf C}_2\langle{\bf B}_2\cdot\rangle_{2})^{-1}=\cI-{\bf C}_2{\bf E}_2^{-1}\langle{\bf
 B}_2\cdot\rangle_{2}.
\]
Multiplying \er{113} by \er{114} we deduce that the operators
$\cA_2$ and $\cI+\cK_1$ are invertible or non-invertible
simultaneously. \BBox

\section{\lb{S2} Example}
\begin{figure}[h]
\center{\includegraphics[width=0.65\linewidth]{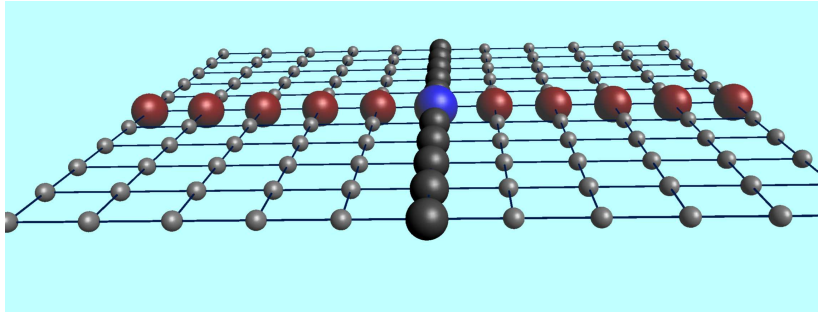}}
\caption{Simple 2D periodic lattice of springs and masses with two
1D defects and with one local defect.} \label{fig2}
\end{figure}

Consider 2D spring-mass lattice with unite masses and unite Hook
modules of springs. We add two perpendicular defects of masses
$1+M_1>0$ and $1+M_2>0$. The mass of cross point of guides is
supposed to be $1+M_1+M_2>0$. Denoting anti-plane displacements at
lattice points ${\bf n}=(x,y)\in\Z^2$ as $u_{\bf n}$,  we can write
the equation of wave motion as
\[\lb{200}
 -(\D_{\rm discr} u)_{\bf n}=\o^2u_{\bf n}+\o^2\ca M_1 u_{\bf n},& if\ y=0,\ x\ne0,\\
                                                   M_2 u_{\bf n},& if\ x=0,\ y\ne0,\\
                                                  (M_1+M_2) u_{\bf
                                                  n},& if\
                                                  x=y=0,
                                               \ac
\]
where $\o$ is a frequency, and the energy $\o^2$ plays the role of
spectral parameter of our discrete periodic (Laplace $\D_{\rm
discr}$) operator with defects. Applying Fourier-Floquet-Bloch
transformation (here it is the Fourier series)
\[\lb{201}
 u({\bf k})=\sum_{{\bf n}\in\Z^2}e^{2\pi i{\bf n}^{\top}{\bf k}}u_{\bf
 n},\ \ {\bf k}=(k_1,k_2)\in[0,1]^2
\]
we rewrite \er{200} as an integral operator
\[\lb{202}
 (4-2\cos2\pi k_1-2\cos2\pi
 k_2)u=\o^2u+\o^2M_2\int_0^1udk_1+\o^2M_1\int_0^1udk_2
\]
or in our notations \er{001}
\[\lb{203}
 Au-\o^2M_2\langle u\rangle_1-\o^2M_1\langle u\rangle_2=0
\]
with $A=4-2\cos2\pi k_1-2\cos2\pi k_2-\o^2$. In fact, the problem
consists of the determining the spectrum (extended eigenvalue
problem for $\o^2$) of the operator from \er{203}. The operator
belongs to $\mA_{2,1}$ and we can use our results for this problem.
Following Definition \ref{D2} introduce (we do not use bold fonts
for scalars)
\[\lb{204}
 C_0=A^{-1},\ \ E_0=A,\ \ C_1=-\o^2M_2A^{-1},\ \ C_2=-\o^2M_1A^{-1},
\]
\[\lb{205}
 E_1(k_2)=1-M_2\langle A^{-1}\rangle_1=1+\o^2M_2\ca \frac{-1}{\sqrt{(2\cos 2\pi k_2-4+\o^2)^2-4}},& if\ \o^2<2-2\cos 2\pi k_2,\\
                                                           \frac{1}{\sqrt{(2\cos 2\pi k_2-4+\o^2)^2-4}},& if\
                                                           \o^2>6-2\cos
                                                           2\pi k_2.
                                                       \ac
\]
If we take $M_1,k_1$ instead of $M_2,k_2$ in \er{205} then we obtain
the identity for $E_2(k_1)$. Using the results of Theorem \ref{T1}
we can describe the spectrum. The continuous part of the spectrum
consists of three components. The first component $\s_0$ \er{013}
corresponds to the propagative waves without attenuation. So, the
energy interval for such waves is
\[\lb{206}
 \s_0=\{\o^2:\ A=0\ for\ some\ {\bf k}\in[0,1]^2\}=[0,8].
\]
The second component $\s_1$ \er{014} corresponds to the guided waves
which propagate along the defect of masses $1+M_2$ and exponentially
decay in perpendicular directions. The energy interval for such
waves is
\[\lb{207}
 \s_1=\{\o^2:\ E_1=0\ for\ some\ k_2\in[0,1]\}=\ca [\frac{4}{1-M_2^2},\frac{6+2\sqrt{8M_2^2+1}}{1-M_2^2}],& M_2<0, \\
                                                  [0,\frac{-6+2\sqrt{8M_2^2+1}}{1-M_2^2}],&
                                                  M_2>0.
 \ac
\]
The third component $\s_2$ \er{014} has the same form as $\s_1$ but
with $M_1$ instead of $M_2$. It is the energy interval for the
guided waves that propagates along the defect of masses $1+M_1$. All
these continuous spectral components are the same as for the lattice
with single line defects, see \cite{K1}. The new thing in our
example is that the crossing of line defects can create the discrete
spectrum. The discrete spectral component is
\[\lb{208}
 \s_3=\{\o^2:\ \cI+\cK_1\ is\ non-invertible\},
\]
where (see \er{009},\er{010},\er{011}, and \er{015})
\[\lb{209}
 \cK_1u=\int_{[0,1]^2}(D_2-D_1)({\bf k},{\bf k}')u({\bf k}')d{\bf k}',\ \ u\in
 L^2_{2,1}\ \ and
\]
\[\lb{210}
 D_1({\bf k},{\bf k}')=C_1({\bf k})E_1^{-1}(k_2)B_1(k_1',k_2)C_2(k_1',k_2)
 B_2({\bf k}'),
\]
\[\lb{211}
 D_2({\bf k},{\bf k}')=C_2({\bf
 k})E_2^{-1}(k_1)\int_0^1B_2(k_1,k_2'')D_1(k_1,k_2'',{\bf
 k}')dk_2''.
\]
Due to non-trivial kernels $D_1,D_2$ the problem of presence or
absence of eigenvalues can be complex and lengthy. Nevertheless,
there are methods that allow to solve this problem effectively.

\section{\lb{S3} Conclusion}
In the current paper we extend some results from \cite{Kjmaa},
\cite{K1a} about the algebra of discrete periodic operators with
parallel defects to the algebra of discrete periodic operators with
perpendicular defects. We did this in the 2D case only. Even in 2D
case we lost the explicit algorithm of finding discrete spectrum.
Now it is not based on simple algebraic operations on some
matrix-valued functions and few number of integrations as it was for
parallel defects. The same thing is expected for multidimensional
periodic operators with various crossing defects. While the
situation is more or less clear in general (abstractly), the
explicit algorithms of finding spectra corresponding to the crossing
defects of lower dimensions are probably not exist (not so simple).

\section*{Acknowledgements}
This work was partially supported by the RSF project
N\textsuperscript{\underline{o}}15-11-30007.

\bibliography{bibl_perp}

\end{document}